\documentclass[conference]{IEEEtran}

\usepackage{cite}
\usepackage{amsmath, mathtools}
\usepackage{amsfonts}
\usepackage{amssymb}
\usepackage{textcomp}

\usepackage[caption=false,font=footnotesize]{subfig}

\usepackage{fixltx2e}

\usepackage{url}
\usepackage{comment}

\usepackage{siunitx}
\usepackage{multirow}
\usepackage{booktabs}
\usepackage{tabularx}

\usepackage[shortlabels]{enumitem}

\usepackage{tikz}
\usetikzlibrary{arrows,decorations.markings,positioning}
\usepackage[european]{circuitikz}
\usepackage{pgfplots}
\usepgfplotslibrary{groupplots}

\pgfplotsset{compat=1.10}
\newlength\figureheight
\newlength\figurewidth

\usepackage{graphicx}
\usetikzlibrary{external}
\tikzexternaldisable

\usepackage{xspace}
\usepackage{accents}
\usepackage{ifthen}
\usepackage{blindtext}

\usepackage{bm}

\newcommand{\tabstretch}{1.45}

\DeclareMathOperator{\real}{Re}

\newcommand{\abs}[1]{{\lvert#1\rvert}}

\newcommand{\mthcal}[1]{{\mathcal{#1}}}

\newcommand{\conj}{^\ast}

\newcommand{\aSrc}[1]{\expandafter\hat#1}
\newcommand{\aDst}[1]{\expandafter\check#1}

\newcommand{\aLB}[1]{\expandafter\underaccent{\bar}#1}
\newcommand{\aUB}[1]{\expandafter\bar#1}

\newcommand{\aAC}[1]{\expandafter\tilde#1}
\newcommand{\aDC}[1]{\expandafter\bar#1}

\newcommand{\aSh}[1]{\expandafter\tilde#1}
\newcommand{\aBar}[1]{\expandafter\bar#1}
\newcommand{\aUBar}[1]{\expandafter\underaccent{\bar}#1}

\newcommand{\aOrg}[1]{#1^{\mathrm{org}}}
\newcommand{\aRed}[1]{#1^{\mathrm{red}}}

\newcommand{\aUndirected}[1]{\expandafter\underaccent{\bar}#1}

\newcommand{\sE}{{\mthcal{E}}}

\newcommand{\sI}{{\mthcal{I}}}

\newcommand{\sV}{{\mthcal{V}}}

\newcommand{\eSrc}{{\aSrc{\epsilon}}}
\newcommand{\eDst}{{\aDst{\epsilon}}}

\newcommand{\hynet}{\emph{hynet}\xspace}

\newcommand{\errinj}{\varepsilon_\mathrm{disp}}
\newcommand{\errflow}{\varepsilon_\mathrm{flow}}

\tikzset{circledletter/.style={circle,draw,inner sep=0mm,minimum size=2.75mm,font=\scriptsize}}

\begin{document}

\title{Feature- and Structure-Preserving Network Reduction for Large-Scale Transmission Grids}

\author{\IEEEauthorblockN{Julia Sistermanns, Matthias Hotz, and Wolfgang Utschick}
\IEEEauthorblockA{Professur f\"ur Methoden der Signalverarbeitung\\
Technische Universit\"at M\"unchen\\
Munich, Germany}%
\and
\IEEEauthorblockN{Dominic Hewes and Rolf Witzmann}
\IEEEauthorblockA{Professur f\"ur Elektrische Energieversorgungsnetze\\
Technische Universit\"at M\"unchen\\
Munich, Germany}%
}

\maketitle

\begin{abstract}
Many countries are currently challenged with the extensive integration of renewable energy sources, which necessitates vast capacity expansion measures. These measures in turn require comprehensive power flow studies, which are often computationally highly demanding. In this work a reduction strategy for large-scale grid models is introduced which not only reduces the model complexity but also preserves the structure and designated grid features. The objective is to ensure that areas crucial to the behavior and the relation of all elements to their physical counterparts remain unchanged. This is accomplished through a specifically designed reduction method for suitable areas identified through topological, electrical and market-based approaches for which we provide an open-source implementation. We show that the proposed strategy adapts to various models and accomplishes a strong reduction of buses and branches while retaining a low dispatch and branch flow deviation. Furthermore, the accuracy of the reduction generalizes well to other scenarios. 
\end{abstract}


\section{Introduction}
\label{sec:introduction}

In several countries, power system operators face the challenges of an extensive integration of renewable energy sources~(RES). This involves a decentralization and geographical dispersion of energy production that generally necessitates capacity expansion measures in the transmission grid~\cite{Capros:2013aa,Flechtner:2015aa}. The decision making for capacity expansion measures relies on comprehensive power flow studies under projected scenarios~\cite{Capros:2013aa,Netzentwicklungsplan2017a}. For large-scale transmission grids such power flow studies involve substantial and potentially prohibitive computational efforts, which require the utilization of network reduction techniques to reduce model complexity while retaining adequate accuracy~\cite{Papaemmanouil:2011aa}. The existing (static) network reduction techniques may be categorized as \emph{electrical equivalencing} and \emph{market-based reduction}.

\emph{Electrical equivalencing methods} divide the system into two subsystems, the internal and external subsystem, where the former remains unmodified while the latter is reduced to a small number of buses (the \emph{essential} buses) that replicate the electrical behavior of the external subsystem~\cite{Liacco:1978aa,Bergen:2000aa,Savulescu:2002aa,Papaemmanouil:2011aa}. A widely employed method for electrical equivalencing is the \emph{Ward equivalent}~\cite{Ward:1949aa} (see also~\cite{Wu:1983aa,Brown:1985aa,Duran:1972aa,Monticelli:1979aa,Deckmann:1980aa}), where all non-essential buses in the external system are reduced via Gaussian elimination on the system of linear equations that relates the nodal voltages and injection currents~\cite{Ward:1949aa,Papaemmanouil:2011aa}. As a consequence of this reduction of the bus admittance matrix, the boundary buses between the internal and external subsystem are interconnected by \emph{artificial} branches in the reduced model. Another popular equivalencing method is the \emph{REI equivalent} introduced by Dimo~\cite{Dimo:1975aa} (see also~\cite{Gavrilas:2008aa,Milano:2009aa,Tinney:1983aa,Housos:1980aa,Papaemmanouil:2011aa}), where the external subsystem is replaced by artificial \emph{radial equivalent independent} (REI) nodes. To this end, related buses in the external system are grouped and replaced by an REI node using the principle of zero power balance, i.e., the injection into the REI node equals the aggregated injection into the respective group of buses. The REI nodes are then connected to the internal subsystem via a radial network of \emph{artificial} branches that replicate the electrical behavior at the considered system state~\cite{Dimo:1975aa,Papaemmanouil:2011aa}.

\emph{Market-based reduction methods} implement zonal aggregations based on results of the optimal power flow (OPF), targeting a reduced model with almost identical power flow among the retained buses. For example, Singh and Srivastava~\cite{Singh:2005aa} proposed an approach that aggregates clusters of buses with almost identical \emph{local marginal prices} (LMPs), as the latter suggests that their aggregation may have a negligible impact on the power flow~\cite{Singh:2005aa,Papaemmanouil:2011aa}. While the LMPs are the means to select a subsystem, the reduction thereof utilizes the REI equivalent. Another approach based on \emph{power transfer distribution factors} (PTDFs) was proposed by Shi and Tylavsky~\cite{Shi:2015aa}, where buses with a similar contribution to designated interzonal power flows are grouped in zones. Each of these zones are aggregated to a single bus and, subsequently, the interzonal flows are modeled using \emph{artificial} branches between zones.

In the context of capacity expansion planning, the introduction of \emph{artificial} branches and buses in the reduced model by the aforementioned methods constitutes a significant drawback. For example, consider that a certain branch in the reduced model is identified, whose capacity should be uprated by adding an additional circuit to that transmission corridor. If the branch is artificial, there is no relation to a physical corridor and the expansion measure cannot be directly realized within the real system. Furthermore, these methods do not explicitly separate the different voltage levels, which potentially leads to the physically unreasonable aggregation of different voltage levels that may incentivize invalid candidate expansion measures. To address these issues, this work proposes a network reduction approach that preserves \emph{features} as well as the \emph{structure} of the grid. By preserving features, like transformers, physically invalid or application-adverse reduction measures are avoided, while preserving the structure ensures that every element in the reduced model possesses a physical counterpart.

In the following, Section~\ref{sec:model} defines the system model and reduction accuracy measure. Section~\ref{sec:features} introduces the notion of features and Section~\ref{sec:reduction} presents the reduction method. Its application is illustrated in Section~\ref{sec:application}, while its robustness with respect to different scenarios and grid models is shown in Section~\ref{sec:verification}. Section~\ref{sec:conclusion} concludes the paper.
An open-source implementation of the proposed reduction method is provided as a toolbox in \hynet, an OPF framework for hybrid AC/DC power systems, which is available at~\cite{Hotz:2019aa}.

\section{System Model and Reduction Accuracy}
\label{sec:model}

This work is based on the system model of \hynet~\cite{Hotz:2018ab},~\cite{Hotz:2019aa} which includes a bus model with shunt compensation, a universal branch model (AC and DC lines and cables, transformers, phase shifters), and a converter model (inverters, rectifiers, voltage source converters). For the sake of simplicity, only those model details that are essential to this work are introduced here. The set of buses is $\sV$, the set of branches is $\sE$, the set of generators\footnote{In \hynet~\cite{Hotz:2018ab},~\cite{Hotz:2019aa} the set of injectors, also denoted by $\sI$, includes conventional and RES-based generation as well as dispatchable and fixed loads. Due to the scope of this work, we focus on the subset of generators.} is $\sI$, and the series impedance of branch $k\in\sE$ is $\aBar{z_k}=1/\aBar{y_k}$. For a given OPF solution, the active power dispatch of generator $i\in\sI$ is $P_i$ and the active power flow on branch $k\in\sE$ is $p_k$, where the latter is considered as the maximum of the absolute active power flow at its two terminals.\footnote{With \hynet's branch model in~\cite[Fig.~1b]{Hotz:2018ab}, the active power flow $p_k$ as considered here is given by $p_k = \max(\abs{\real(V_{\eSrc(k)}\aSrc{I_k}\conj)}, \abs{\real(V_{\eDst(k)}\aDst{I_k}\conj)})$.} Furthermore, the dual variable of the active power balance equation of bus $n\in\sV$ is $\lambda_n$, which corresponds to the LMP in case of a zero duality gap of the OPF problem~\cite{Hotz:2018aa}.

According to the application in capacity expansion planning, the \emph{accuracy} of the network reduction is considered in terms of the similarity of the \emph{generator dispatch} and \emph{branch flows} produced by an OPF calculation. To this end, let $\aOrg{P_i}$ and $\aRed{P_i}$ be the active power dispatch of generator $i\in\sI$ and $\aOrg{p_k}$ and $\aRed{p_k}$ the active power flow on branch $k\in\aRed{\sE}\subset\aOrg{\sE}$ in the original and reduced model, respectively, where $\aRed{\sE}$ is the set of branches in the reduced model. Therewith, the \emph{contribution-weighted} \emph{mean relative dispatch error} $\errinj$ and \emph{mean relative flow error} $\errflow$ reads
\vspace{0.2cm}
\begin{align}
    \errinj &= \frac{1}{\sum_{i\in\sI}\abs{\aOrg{P_i}}} \sum_{i\in\sI} \abs{\aRed{P_i} - \aOrg{P_i}}\\
    \errflow &= \frac{1}{\sum_{k\in\aRed{\sE}}\aOrg{p_k}} \sum_{k\in\aRed{\sE}} \abs{\aRed{p_k} - \aOrg{p_k}}\,.
\end{align}

\section{The Notion of Features}
\label{sec:features}

\emph{Features} are defined as entities in the model which are essential to the application-relevant accuracy and validity of the derived results and conclusions. With regard to capacity expansion planning, the following \emph{features} can be identified:
\begin{enumerate}[(a)]
    \item \emph{Transformers}, as their reduction may result in the aggregation of different voltage levels.
    \item \emph{Converters}, as their reduction may result in the aggregation of AC and DC grids.
    \item Selected \emph{branches} with a particular relevance to capacity expansion decisions, namely
        \begin{itemize}
            \item highly loaded and congested branches, i.e., their flow is observed to be close to or at the capacity rating, and
            \item branches of long transmission lines, e.g., longer than $50\,$km, as they are typically associated with bulk transmission.
        \end{itemize}
    \item \emph{Terminal buses of conventional generators}, as their power dispatch impacts the power flow significantly. This includes the reference bus.
\end{enumerate}
By preserving such features during the reduction process, application-specific requirements on the reduced model can be incorporated. Furthermore, as illustrated later on, empirically determined features may be introduced to control and improve the accuracy of the model reduction.

\section{Feature- and Structure-Preserving\\ Network Reduction}
\label{sec:reduction}

In order to arrive at a feature- and structure-preserving network reduction, the reduction process must be (a) aware of features and (b) retain the relation of every entity in the reduced model to its physical counterpart. To this end, the proposed method aims at the identification of a multitude of small \emph{subgrids} within the transmission grid that are suited for reduction. These subgrids are then filtered based on features, i.e., if a subgrid contains a feature it is excluded from the reduction process. Subsequently, the remaining subgrids are aggregated in a fashion that avoids the introduction of artificial entities. The following description of the proposed method starts with the subgrid reduction and, subsequently, introduces three approaches to select suitable subgrids.

\subsection{Subgrid Reduction}
\label{subgrid_red}
As the avoidance of artificial entities in the reduced model is considered essential, it is assumed that the subgrid can be reduced to a \emph{single} bus, where the latter is the \emph{representative bus} selected among the buses of the subgrid. Compared to the Ward and REI equivalent, this is a very restricted reduction process in terms of modeling the subgrid's impact on the surrounding grid. In the proposed method, this limitation is compensated by the selection of comparably \emph{small} subgrids, whose impact on the electrical behavior of the overall system is negligible.

In case that the considered subgrid does \emph{not} contain any features, the following steps are performed to reduce the subgrid to its representative bus:
\begin{enumerate}[(a)]
    \item The terminal bus of any generators and loads within the subgrid is set to the representative bus.
    \item All reactive power compensation (shunts) within the subgrid is moved to the representative bus.
    \item All branches or converters that connect a bus of the surrounding grid to a bus within the subgrid are connected to the representative bus instead.
    \item All remaining buses and branches of the subgrid are removed. The line charging of the removed branches is modeled as a shunt at the representative bus.
\end{enumerate}
Due to this comparably coarse subgrid reduction, a careful selection of subgrids is essential. To this end, we utilize insights into the system behavior to identify subgrids with a potentially negligible impact on the electrical behavior of the overall system. Hereafter, three approaches are presented, which are based on topological, electrical, and market insights.

\subsection{Topology-Based Subgrid Selection}

Transmission grids often exhibit small subgrids at the boundary of the grid which are only connected by a single corridor, i.e., one branch or several parallel branches. These subgrids include single buses, lines of buses and small ``islands'', i.e., small groups of buses which are connected to the main grid via a single corridor. Such structures may contain several branches and be meshed, but often their internal power flow is not crucial to the overall grid behavior.

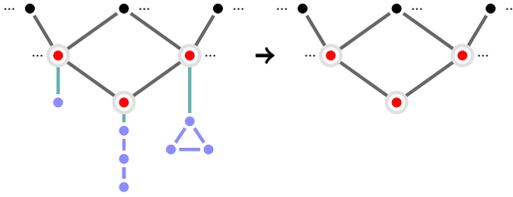
\begin{figure}[hb]
\centering
\tikzsetnextfilename{topology_reduction}
\begin{tikzpicture}[scale=0.25, transform shape, line width=0.45mm,
	bus/.style={circle,inner sep=2.5mm,draw=white,fill=black},
	connecting/.style={circle,inner sep=2.5mm,draw=white,fill=red},
	island/.style={circle,inner sep=2.5mm,draw=white,fill=blue!45!white},
	connector/.style={draw=teal!60!white},
	islandlink/.style={draw=blue!45!white},
	aclink/.style={draw=black!60!white},
	dclink/.style={dash pattern=on 4.625pt off 2.375pt,draw=black!50!white}] 
%
%
\node (bus01) at ( 0, 0) {};
\node (bus02) at (-3.5, 2.5) {};\node[left=5mm of bus02] {\Huge {...}};
\node (bus03) at (3.5, 2.5) {};\node[right=5mm of bus03] {\Huge {...}};
\node (bus04) at ( 0, 5) {};\node[right=5mm of bus04] {\Huge {...}};
\node (bus13) at (-5, 5) {};\node[left=5mm of bus13] {\Huge {...}};
\node (bus14) at (5, 5) {};\node[right=5mm of bus14] {\Huge {...}};
\node (bus05) at ( 0, -1.5) {};
\node (bus10) at ( 0, -3) {};
\node (bus11) at ( 0, -4.5) {};
\node (bus12) at (-3.5, 0) {};
\node (bus17) at (3.5, -1) {};
\node (bus18) at (2.5, -2.5) {};
\node (bus19) at (4.5, -2.5) {};
\node (branch01) at (0, -2.5){}; 
\node (bus06) at (14.5, 0) {};
\node (bus07) at (11, 2.5) {};\node[left=5mm of bus07] {\Huge {...}};
\node (bus08) at (18, 2.5) {};\node[right=5mm of bus08] {\Huge {...}};
\node (bus09) at ( 14.5, 5) {};\node[right=5mm of bus09] {\Huge {...}};
\node (bus15) at (9.5, 5) {};\node[left=5mm of bus15] {\Huge {...}};
\node (bus16) at (19.5, 5) {};\node[right=5mm of bus16] {\Huge {...}};
%
%
\draw[aclink] (bus01) -- (bus02) -- (bus04) -- (bus03)--(bus01);
\draw[aclink] (bus02) -- (bus13);
\draw[aclink] (bus03) -- (bus14);
\draw[connector] (bus01) -- (bus05);
\draw[connector] (bus02) -- (bus12);
\draw[connector] (bus03) -- (bus17);
\draw [->](7,2.5) -- (8,2.5);
\draw[islandlink] (bus05) -- (bus10) -- (bus11) -- (bus10);
\draw[aclink] (bus06) -- (bus07) -- (bus09) -- (bus08)--(bus06);
\draw[aclink] (bus07) -- (bus15);
\draw[aclink] (bus08) -- (bus16);
\draw[islandlink] (bus17) -- (bus18) -- (bus19) -- (bus17);
%
%
\node[connecting, inner sep=3.8mm, draw=black!12, fill=black!12] at (bus01) {};
\node[connecting] at (bus01) {};
\node[connecting, inner sep=3.8mm, draw=black!12, fill=black!12] at (bus02) {};
\node[connecting] at (bus02) {};
\node[connecting, inner sep=3.8mm, draw=black!12, fill=black!12] at (bus03) {};
\node[connecting] at (bus03) {};
\node[bus] at (bus04) {};
\node[bus] at (bus13) {};
\node[bus] at (bus14) {};
\node[island] at (bus05) {};
\node[island] at (bus17) {};
\node[island] at (bus18) {};
\node[island] at (bus19) {};
\node[connecting, inner sep=3.8mm, draw=black!12, fill=black!12] at (bus06) {};
\node[connecting] at (bus06) {};
\node[connecting, inner sep=3.8mm, draw=black!12, fill=black!12] at (bus07) {};
\node[connecting] at (bus07) {};
\node[connecting, inner sep=3.8mm, draw=black!12, fill=black!12] at (bus08) {};
\node[connecting] at (bus08) {};
\node[bus] at (bus09) {};
\node[bus] at (bus15) {};
\node[bus] at (bus16) {};
\node[island] at (bus10) {};
\node[island] at (bus11) {};
\node[island] at (bus12) {};
%
%
\end{tikzpicture}
\caption{Reduction of single buses, lines of buses and small ``islands'' at the boundary of the grid.}
\label{topology_reduction}
\vspace{-0.1cm}
\end{figure}

Consequently, such a subgrid is a suitable candidate for reduction, using the shared bus with the main grid as the subgrid's representative bus as illustrated in Fig.~\ref{topology_reduction}. 

\subsection{Electrical Coupling-Based Subgrid Selection}

When two buses are connected via a branch with a very low series impedance, their electrical states are strongly coupled compared to buses linked by medium to high impedance branches. Therefore their aggregation may not affect the overall system behavior significantly and subsequently they are potentially suited for reduction.

\begin{figure}[t]
\centering
\tikzsetnextfilename{electric_reduction}
\begin{tikzpicture}[scale=0.25, transform shape, line width=0.45mm,
	bus/.style={circle,inner sep=2.5mm,draw=white,fill=black},
	connecting/.style={circle,inner sep=2.5mm,draw=white,fill=red},
	island/.style={circle,inner sep=2.5mm,draw=white,fill=blue!45!white},
	connector/.style={draw=teal!60!white},
	islandlink/.style={draw=blue!45!white},
	aclink/.style={draw=black!60!white},
	dclink/.style={dash pattern=on 4.625pt off 2.375pt,draw=black!50!white}] 

 Place buses

\node (bus01) at ( 0, 0) {};\node[left=5mm of bus01] {\Huge {...}};
\node (bus02) at (0, -2.5) {};\node[left=5mm of bus02] {\Huge {...}};
\node (bus03) at (0, -5) {};\node[left=5mm of bus03] {\Huge {...}};
\node (bus04) at ( 2, -2.5) {};
\node (bus05) at (6, -2.5) {};
\node (bus06) at (8, 0) {};\node[right=5mm of bus06] {\Huge {...}};
\node (bus07) at ( 8, -2.5) {};\node[right=5mm of bus07] {\Huge {...}};
\node (bus08) at ( 8, -5) {};\node[right=5mm of bus08] {\Huge {...}};
\node (b01) at (14, 0) {};\node[left=5mm of b01] {\Huge {...}};
\node (b02) at (14, -2.5) {};\node[left=5mm of b02] {\Huge {...}};
\node (b03) at (14, -5) {};\node[left=5mm of b03] {\Huge {...}};
\node (b04) at (16, -2.5) {};
\node (b06) at (18, 0) {};\node[right=5mm of b06] {\Huge {...}};
\node (b07) at (18, -2.5) {};\node[right=5mm of b07] {\Huge {...}};
\node (b08) at (18, -5) {};\node[right=5mm of b08] {\Huge {...}};

 AC links

\draw[aclink] (bus01) -- (bus04);
\draw[aclink] (bus02) -- (bus04);
\draw[aclink] (bus03) -- (bus04);
\draw[connector] (bus04) -- (bus05);
\draw[aclink] (bus06) -- (bus05);
\draw[aclink] (bus07) -- (bus05);
\draw[aclink] (bus08) -- (bus05);
\draw [->](10.5, -2.5) -- (11.5, -2.5);
\draw[aclink] (b01) -- (b04);
\draw[aclink] (b02) -- (b04);
\draw[aclink] (b03) -- (b04);
\draw[aclink] (b06) -- (b04);
\draw[aclink] (b07) -- (b04);
\draw[aclink] (b08) -- (b04);
%
%
\node[connecting, inner sep=3.8mm, draw=black!12, fill=black!12] at (bus04) {};
\node[connecting] at (bus04) {};
\node[connecting, inner sep=3.8mm, draw=black!12, fill=black!12] at (b04) {};
\node[connecting] at (b04) {};
\node[bus] at (bus01) {};
\node[bus] at (bus02) {};
\node[bus] at (bus03) {};
\node[island] at (bus05) {};
\node[bus] at (bus06) {};
\node[bus] at (bus07) {};
\node[bus] at (bus08) {};
\node[bus] at (b01) {};
\node[bus] at (b02) {};
\node[bus] at (b03) {};
\node[bus] at (b06) {};
\node[bus] at (b07) {};
\node[bus] at (b08) {};

%
%
\end{tikzpicture}
\caption{Reduction of branches with a very low series impedance.}
\label{electric_reduction}
\vspace{-0.2cm}
\end{figure}
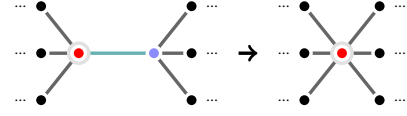

Such branches $k\in\sE$ are identified by comparing their series admittance in Ohms to a threshold parametrized by $\tau \in [0,1]$ which is relative to the maximum series impedance in the grid:\footnote{In case of parallel branches, the equivalent series impedance must be considered. For the sake of a simple presentation this is not elaborated here.}
\begin{equation}
\centering
|\aBar{z_k}| \leq \tau \cdot \max\limits_{k' \in \sE}|\aBar{z_{k'}}|\,.
\end{equation}
The subgrid associated with such a  branch consists of the branch itself as well as its terminal buses. It is reduced to its representative bus, which is either one of the terminals, see Fig.~\ref{electric_reduction}.\\

\subsection{Market-Based Subgrid Selection}

Similar to the market-based reduction of Singh and Srivastava~\cite{Singh:2005aa} the following approach utilizes LMPs to identify candidate subgrids. However, here the focus is on other insights gained from similar LMPs, namely the following two valuable implications: Firstly, a connected group of buses with almost identical LMPs indicates an area with low electrical losses. Secondly, if congestion occurs in the grid, the LMPs in the vicinity of the congested branches usually diverge. As the overall system behavior is sensitive to areas with high losses and congestion, these should be preserved. On the contrary, the internal power flow in areas with similar LMPs is potentially negligible, rendering them candidates for reduction.

Such subgrids are identified through a breadth first search approach starting from candidate buses with two or more connected corridors. In this clustering process, the reference is always the LMP of the candidate bus where the search started. If a cluster is identified, this candidate will be chosen as the respective representative bus. A bus $n$ that is connected to a representative bus $r$ is included in the respective subgrid if the deviation of their LMPs\footnote{Note that the dual variables $\lambda_n$ only equal the LMPs if the associated primal solution is globally optimal and if the OPF exhibits a zero duality gap, cf. e.g.~\cite{Hotz:2018aa}. However, with respect to the subgrid identification, this issue is of minor relevance and, thus, not considered to simplify the discussion.} is below the threshold $\delta>0$, i.e.,
\begin{align}
   \abs{\lambda_n - \lambda_r} \leq \delta\,.
\end{align}
The process is repeated iteratively for all buses that are connected to a subgrid until all boundary buses of the subgrids connect to buses of the main grid that exhibit an LMP deviation beyond $\delta$ to the respective representative bus. After the clustering process, any overlapping subgrids are combined, using one of their representative buses as the representative bus of the union. All subgrids are then subject to the reduction process in Section~\ref{subgrid_red}.
\vspace{-0.2cm}
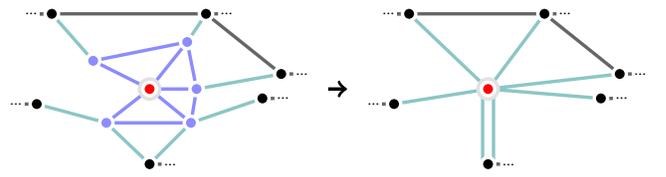
\begin{figure}[hb!]
\centering
\tikzsetnextfilename{clusters}
\begin{tikzpicture}[scale=0.25, transform shape, line width=0.45mm,
	bus/.style={circle,inner sep=2.5mm,draw=white,fill=black},
	connecting/.style={circle,inner sep=2.5mm,draw=white,fill=red},
	island/.style={circle,inner sep=2.5mm,draw=white,fill=black},
	d1/.style={circle,inner sep=2.5mm,draw=white,fill=blue!45!white},
	connector/.style={draw=teal!60!white},
	islandlink/.style={draw=blue!45!white},
	greenlink/.style={draw=teal!45!white},
	aclink/.style={draw=black!60!white},
	dclink/.style={dash pattern=on 4.625pt off 2.375pt,draw=black!50!white}] 
%
\node (bus01) at (0, 0) {};
\node (bus02) at (2.5, 0) {};
\node (bus03) at ( 2.2, -1.8) {};
\node (bus04) at ( -2.3, -1.8) {};
\node (bus05) at ( 0,-4) {};
\node (bus06) at ( -3, 1.5) {};
\node (bus07) at ( 2, 2.5) {};
\node (bus08) at ( 3, 4) {};
\node (bus09) at ( 6, -0.5) {};
\node (bus10) at ( -6, -0.8) {};
\node (bus11) at ( -5.2, 4) {};
\node (bus12) at (7, 0.8) {};

\node (b01) at (18, 0) {};
\node (b05) at ( 18,-4) {};
\node (b08) at ( 21, 4) {};
\node (b09) at ( 24, -0.5) {};
\node (b10) at ( 13, -0.8) {};
\node (b11) at ( 13.8, 4) {};
\node (b12) at (25, 0.8) {};

\node [left=5mm of bus11] (r1) {\Huge {...}};
\node [left=5mm of bus10] (r2) {\Huge {...}};
\node [right=5mm of bus08] (r3) {\Huge {...}};
\node [right=5mm of bus09] (r4) {\Huge {...}};
\node [right=5mm of bus05] (r5) {\Huge {...}};
\node [right=5mm of bus12] (r6) {\Huge {...}};

\node [left=5mm of b11] (rr1) {\Huge {...}};
\node [left=5mm of b10] (rr2) {\Huge {...}};
\node [right=5mm of b08] (rr3) {\Huge {...}};
\node [right=5mm of b09] (rr4) {\Huge {...}};
\node [right=5mm of b05] (rr5) {\Huge {...}};
\node [right=5mm of b12] (rr6) {\Huge {...}};

\node [right=0.1mm of b01] (n1) {};
\node [right=0.1mm of b05] (n5) {};
\node [left=0.1mm of b01] (l1) {};
\node [left=0.1mm of b05] (l5) {};
%
%
\draw[aclink] (bus11) -- (r1);
\draw[aclink] (bus10) -- (r2);
\draw[aclink] (bus08) -- (r3);
\draw[aclink] (bus09) -- (r4);
\draw[aclink] (bus05) -- (r5);
\draw[aclink] (bus11) -- (r1);
\draw[aclink] (bus12) -- (r6);
\draw[islandlink] (bus01) -- (bus02);
\draw[islandlink] (bus01) -- (bus03);
\draw[islandlink] (bus01) -- (bus04);
\draw[islandlink] (bus01) -- (bus06);
\draw[islandlink] (bus01) -- (bus07);
\draw[islandlink] (bus02) -- (bus07);
\draw[islandlink] (bus03) -- (bus04);
\draw[islandlink] (bus03) -- (bus02);
\draw[islandlink] (bus06) -- (bus07);
\draw[greenlink] (bus02) -- (bus12); 
\draw[greenlink] (bus09) -- (bus03) -- (bus05) -- (bus04) -- (bus10);
\draw[greenlink] (bus06) -- (bus11); 
\draw[greenlink] (bus08) -- (bus07);
\draw[aclink] (bus11) --(bus08) -- (bus12);

\draw [->](9.5,0) -- (10.5,0);

\draw[aclink] (b11) -- (rr1);
\draw[aclink] (b10) -- (rr2);
\draw[aclink] (b08) -- (rr3);
\draw[aclink] (b09) -- (rr4);
\draw[aclink] (b05) -- (rr5);
\draw[aclink] (b11) -- (rr1);
\draw[aclink] (b12) -- (rr6);
\draw[greenlink] (b01) -- (b09);
\draw[greenlink] (b01) -- (b10);
\draw[greenlink] (b01) -- (b11);
\draw[greenlink] (b01) -- (b12);
\draw[greenlink] (b01) -- (b08);
\draw[greenlink] (l1) -- (l5);
\draw[greenlink] (n1) -- (n5);
\draw[aclink] (b11) -- (b08) -- (b12);

%
\node[connecting, inner sep=3.8mm, draw=black!12, fill=black!12] at (bus01) {};
\node[connecting] at (bus01) {};
\node[d1] at (bus02) {};
\node[d1] at (bus03) {};
\node[d1] at (bus04) {};
\node[island] at (bus05) {};
\node[d1] at (bus06) {};
\node[d1] at (bus07) {};
\node[island] at (bus08) {};
\node[island] at (bus09) {};
\node[island] at (bus10) {};
\node[island] at (bus11) {};
\node[island] at (bus12) {};
\node[connecting, inner sep=3.8mm, draw=black!12, fill=black!12] at (b01) {};
\node[connecting] at (b01) {};
\node[island] at (b05) {};
\node[island] at (b08) {};
\node[island] at (b09) {};
\node[island] at (b10) {};
\node[island] at (b11) {};
\node[island] at (b12) {};

%
%
\end{tikzpicture}
\caption{Reduction of a connected group of buses with similar LMPs.}
\label{b2.1.7}
\end{figure}

\section{Reduction Parameter Selection}
\label{sec:application}

The proposed reduction method depends on a set of tuning parameters. To better illustrate their impact on the result and to show how to arrive at a proper parametrization for a desired accuracy, the proposed method is demonstrated in detail for an exemplary grid.
This example can also be found as a tutorial in the open-source implementation provided in~\cite{Hotz:2019aa}.
To this end, the German high voltage transmission grid for the year 2030 as proposed in the network development plan~\cite{Hewes:2016aa},~\cite{Netzentwicklungsplan2017a} is considered. The model consists of a total of 1524 buses and 2208 branches. With the criteria defined in Section~\ref{sec:features}, 856 features are identified.
\vspace{-0.2cm}
\begin{figure}[!t]
\centering
\setlength\figurewidth{0.9\columnwidth}
\tikzsetnextfilename{elec_DE}
\begin{tikzpicture}

\begin{groupplot}[group style={vertical sep=20pt, group size=1 by 2}]
\nextgroupplot[
width=\figurewidth,
height=\figureheight,
label style = {font=\footnotesize}, 
ylabel={Reduction in \%},
xmin=0, xmax=0.1,
ymin=0, ymax=65,
ticklabel style = {font=\footnotesize},
xtick={0, 0.025, 0.05, 0.075, 0.1},
xticklabels={$0$, $0.025$, $0.05$, $0.075$, $0.1$},  
tick align=outside,
tick pos=left,
x grid style={white!69.01960784313725!black},
y grid style={white!69.01960784313725!black}
]

\addplot [black, mark=*, mark size=1, mark options={solid}]
table {%
0 13.3858267716535
0.005 19.5538057742782
0.01 24.2782152230971
0.015 27.755905511811
0.02 32.6115485564304
0.025 35.4330708661417
0.03 38.9107611548556
0.035 42.5853018372703
0.04 44.9475065616798
0.045 46.8503937007874
0.05 49.0157480314961
0.06 53.1496062992126
0.07 56.5616797900262
0.08 59.4488188976378
0.09 62.3359580052493
0.1 64.501312335958
};
\addplot [white!50!black, mark=*, mark size=1, mark options={solid}]
table {%
0 10.8242753623188
0.005 15.8967391304348
0.01 20.3351449275362
0.015 23.4601449275362
0.02 27.6268115942029
0.025 30.2536231884058
0.03 33.3333333333333
0.035 36.5489130434783
0.04 38.5869565217391
0.045 40.3985507246377
0.05 42.2554347826087
0.06 45.9239130434783
0.07 49.2300724637681
0.08 51.9021739130435
0.09 54.0760869565217
0.1 56.0235507246377
};
\addplot [white!50!black, dashed, mark=*, mark size=1, mark options={solid}]
table {%
0 0.21978021978022
0.005 4.61538461538462
0.01 7.69230769230769
0.015 11.4285714285714
0.02 14.9450549450549
0.025 16.9230769230769
0.03 20.4395604395604
0.035 23.5164835164835
0.04 26.8131868131868
0.045 29.010989010989
0.05 30.5494505494506
0.06 33.8461538461538
0.07 38.6813186813187
0.08 43.0769230769231
0.09 47.9120879120879
0.1 50.1098901098901
};
\path [draw=black, fill opacity=0] (axis cs:0,0)
--(axis cs:0,65);
\nextgroupplot[
width=\figurewidth,
height=\figureheight,
label style = {font=\footnotesize}, 
xlabel={Threshold $\tau$},
ylabel={Error in \%},
xmin=0, xmax=0.1,
ymin=0, ymax=30,
ticklabel style = {font=\footnotesize},
xtick={0, 0.025, 0.05, 0.075, 0.1},
xticklabels={$0$, $0.025$, $0.05$, $0.075$, $0.1$},  
tick align=outside,
tick pos=left,
x grid style={white!69.01960784313725!black},
y grid style={white!69.01960784313725!black}
]
\addplot [black, mark=*, mark size=1, mark options={solid}]
table {%
0 0.118546453899469
0.005 0.752489649621045
0.01 1.90498461673772
0.015 3.826360804334
0.02 3.83800998456309
0.025 4.25159110954843
0.03 4.45771200401316
0.035 6.84348464020489
0.04 7.87105066806011
0.045 7.44920035519624
0.05 8.71225273114225
0.06 8.97974389676063
0.07 13.1681392111098
0.08 12.2653326151649
0.09 10.3308305552334
0.1 14.0480698624572
};
\addplot [white!50!black, mark=*, mark size=1, mark options={solid}]
table {%
0 0.070008897005993
0.005 0.577738265913165
0.01 3.14951227081565
0.015 4.1771105675327
0.02 3.94978291166884
0.025 4.94322003554053
0.03 5.81616044051919
0.035 9.4939016972098
0.04 10.4833707361713
0.045 10.2153954577921
0.05 12.4947782079961
0.06 13.6079511683474
0.07 19.5263762208953
0.08 19.6138924814717
0.09 20.0579183160221
0.1 24.3027560329555
};
\path [draw=black, fill opacity=0] (axis cs:0,0)
--(axis cs:0,20);

\end{groupplot}

\end{tikzpicture}
\caption{Result of the electrical coupling-based reduction for different choices of $\tau$. The upper figure shows the reduction of buses (black), branches (gray), and cycles (dashed), while the lower depicts $\errinj$ (black) and $\errflow$ (gray).}
\label{elec_DE}
\vspace{-0.5cm}
\end{figure}
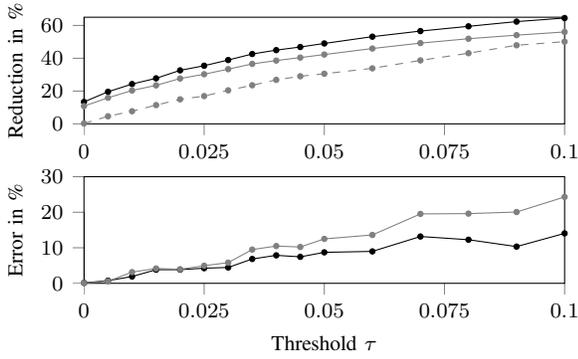

\subsection{Modular Principle of the Reduction Strategy}
\label{combined}
The proposed method provides three different subgrid selection approaches and follows a modular principle: Each approach can be used individually and tuned to the grid model. The parametrization process also offers valuable insights into the grid that can be used to refine the reduction process.

By applying a combined reduction strategy, the degrees of freedom can be used to target a certain accuracy, e.g. a dispatch error $\errinj < 2\%$. We identified the following sequence as the most appropriate ordering:

\begin{itemize}[leftmargin=2.85em]
\item[(i)] Topology-based reduction 
\item[(ii.a)] Electrical coupling-based reduction with threshold $\tau$
\item[(ii.b)] Feature refinement and repetition of (ii.a)
\item[(iii)] Market-based reduction with parameter $\delta$
\end{itemize}

To start with the topology-based reduction is obvious as it does not impact the other reduction steps, is basically independent of parameter adjustments, and introduces only a small error. The ordering of the next steps is not as intuitive, but several experiments have shown that the market-based reduction profits from the electrical model-based reduction and its feature addition, while in the reversed ordering the feature addition is less effective.
\vspace{-0.2cm}

\subsection{Parametrization}
The parametrization of the combined strategy is performed consecutively with each reduction step continuing on the result of the previous step.

The \textit{topology-based reduction} depends on a definition of a ``small'' group of buses. Here 1\% of the total number of buses were considered small. In several studies we found that above a certain threshold, the precise definition of ``small'' has little impact as the boundary structures themselves are small, e.g. for the considered model the choice of 20 or 800 buses leads to the same result. Even though the subgrids are small, the reduction potential is considerable: $13.4 \%$ of all buses and $10.8 \%$ of all branches are reduced by this approach, while the induced error is very small with $\errinj = 0.11\%$ and $\errflow = 0.07\%$.

In contrast, the \textit{electrical coupling-based reduction} depends heavily on its parameter $\tau \in [0,1]$.
When performing a reduction for different $\tau$, it becomes evident that this method offers a trade-off between the reduction and the error: A higher value for $\tau$ enables a more extensive reduction, but usually increases the error for both the dispatch and branch flows, see Fig.~\ref{elec_DE}. 
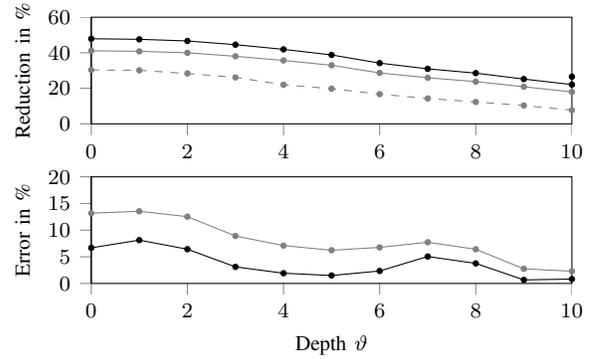
\begin{figure}[!t]
\centering
\tikzsetnextfilename{1shot_DE}
\setlength\figurewidth{0.9\columnwidth}
\begin{tikzpicture}

\begin{groupplot}[group style={vertical sep=20pt, group size=1 by 2}]
\nextgroupplot[
width=\figurewidth,
height=\figureheight,
label style = {font=\footnotesize}, 
ylabel={Reduction in \%},
xmin=0, xmax=10,
ymin=0, ymax=60,
ticklabel style = {font=\footnotesize},
tick align=outside,
tick pos=left,
x grid style={white!69.01960784313725!black},
y grid style={white!69.01960784313725!black}
]

\addplot [black, mark=*, mark size=1, mark options={solid}]
table {%
0 47.9002624671916
1 47.5721784776903
2 46.6535433070866
3 44.5538057742782
4 41.9291338582677
5 38.7795275590551
6 34.1863517060367
7 30.9711286089239
8 28.5433070866142
9 25.1968503937008
10 22.1128608923885
10 26.5748031496063
};
\addplot [white!50!black, mark=*, mark size=1, mark options={solid}]
table {%
0 41.1231884057971
1 40.8061594202899
2 39.9909420289855
3 38.0434782608696
4 35.6884057971014
5 32.9710144927536
6 28.6684782608696
7 25.9510869565217
8 23.7771739130435
9 20.8786231884058
10 17.9800724637681
};
\addplot [white!50!black, dashed, mark=*, mark size=1, mark options={solid}]
table {%
0 30.3296703296703
1 30.1098901098901
2 28.3516483516484
3 26.1538461538462
4 21.978021978022
5 19.7802197802198
6 16.7032967032967
7 14.2857142857143
8 12.3076923076923
9 10.3296703296703
10 7.69230769230769
};
\path [draw=black, fill opacity=0] (axis cs:0,0)
--(axis cs:0,60);

\nextgroupplot[
width=\figurewidth,
height=\figureheight,
label style = {font=\footnotesize}, 
xlabel={Depth $\vartheta$},
ylabel={Error in \%},
xmin=0, xmax=10,
ymin=0, ymax=20,
ticklabel style = {font=\footnotesize},
tick align=outside,
tick pos=left,
x grid style={white!69.01960784313725!black},
y grid style={white!69.01960784313725!black}
]

\addplot [black, mark=*, mark size=1, mark options={solid}]
table {%
0 6.67885090667879
1 8.12627058309718
2 6.42352699080081
3 3.11297744566834
4 1.92266704062884
5 1.49127253815231
6 2.33965935700818
7 5.05542702914279
8 3.76445051772796
9 0.666875090805221
10 0.822759862694402
};
\addplot [white!50!black, mark=*, mark size=1, mark options={solid}]
table {%
0 13.1776636478679
1 13.5380091260137
2 12.5293577735284
3 8.90166353130969
4 7.10874368838231
5 6.22971352854188
6 6.76296924102375
7 7.73296397435898
8 6.41704304352447
9 2.74053444458385
10 2.29485665061789
};
\path [draw=black, fill opacity=0] (axis cs:0,0)
--(axis cs:0,30);

\end{groupplot}

\end{tikzpicture}
\caption{Results of the electrical coupling-based reduction with $\tau = 0.05$ and additional features for different depths $\vartheta$ from critical generators. The upper figure shows the reduction of buses (black), branches (gray), and cycles (dashed), while the lower depicts $\errinj$ (black) and $\errflow$ (gray).}
\label{1shot_DE}
\end{figure}
\begin{figure}[!t]
\centering
\tikzsetnextfilename{cluster_DE}
\setlength\figurewidth{0.9\columnwidth}
\begin{tikzpicture}

\definecolor{color0}{rgb}{0.0156862745098039,0.454901960784314,0.584313725490196}
\definecolor{color1}{rgb}{0.0823529411764706,0.690196078431373,0.101960784313725}

\begin{groupplot}[group style={vertical sep=20pt, group size=1 by 2}]
\nextgroupplot[
width=\figurewidth,
height=\figureheight,
label style = {font=\footnotesize}, 
ylabel={Reduction in \%},
xmin=0, xmax=0.2,
ymin=0, ymax=60,
ticklabel style = {font=\footnotesize},
xtick={0, 0.05, 0.1, 0.15, 0.2},
xticklabels={$0$, $0.05$, $0.1$, $0.15$, $0.2$},
tick align=outside,
tick pos=left,
x grid style={white!69.01960784313725!black},
y grid style={white!69.01960784313725!black}
]
\addplot [black, mark=*, mark size=1, mark options={solid}]
table {%
0 41.9291338582677
0.02 47.4409448818898
0.04 50
0.06 51.1154855643045
0.08 52.1653543307087
0.1 53.3464566929134
0.12 54.002624671916
0.14 54.6587926509186
0.16 55.1181102362205
0.18 55.1837270341207
0.2 55.8398950131234
};
\addplot [white!50!black, mark=*, mark size=1, mark options={solid}]
table {%
0 35.6884057971014
0.02 40.2173913043478
0.04 42.4365942028986
0.06 43.3876811594203
0.08 44.1123188405797
0.1 45.1539855072464
0.12 45.8333333333333
0.14 46.4221014492754
0.16 46.7391304347826
0.18 46.875
0.2 47.3731884057971
};
\addplot [white!50!black, dashed, mark=*, mark size=1, mark options={solid}]
table {%
0 21.978021978022
0.02 28.3516483516484
0.04 32.3076923076923
0.06 34.7252747252747
0.08 35.8241758241758
0.1 38.9010989010989
0.12 40
0.14 41.0989010989011
0.16 41.978021978022
0.18 41.978021978022
0.2 43.2967032967033
};
\path [draw=black, fill opacity=0] (axis cs:0,0)
--(axis cs:0,60);

\nextgroupplot[
width=\figurewidth,
height=\figureheight,
label style = {font=\footnotesize}, 
xlabel={Maximum LMP deviation $\delta$},
ylabel={Error in \%},
xmin=0, xmax=0.2,
ymin=0, ymax=15,
ticklabel style = {font=\footnotesize},
xtick={0, 0.05, 0.1, 0.15, 0.2},
xticklabels={$0$, $0.05$, $0.1$, $0.15$, $0.2$},  
tick align=outside,
tick pos=left,
x grid style={white!69.01960784313725!black},
y grid style={white!69.01960784313725!black}
]

\addplot [black, mark=*, mark size=1, mark options={solid}]
table {%
0 1.92266704062884
0.02 1.9529701934226
0.04 1.52754640110025
0.06 1.51696188145836
0.08 1.51848142740992
0.1 2.07402581298035
0.12 2.18549004697287
0.14 2.62888170816701
0.16 2.95510800114399
0.18 3.36926063583521
0.2 3.40472347486093
};
\addplot [white!50!black, mark=*, mark size=1, mark options={solid}]
table {%
0 7.10874368838231
0.02 8.62141441705108
0.04 9.11030371701055
0.06 9.37522389308101
0.08 9.4062221797721
0.1 10.7662018902574
0.12 11.051832815379
0.14 11.4382489999383
0.16 12.343566460713
0.18 12.9033449722569
0.2 13.3085369784038
};
\path [draw=black, fill opacity=0] (axis cs:0,0)
--(axis cs:0,15);

\end{groupplot}

\end{tikzpicture}
\caption{Results of the market-based reduction for different $\delta$. The upper figure shows the reduction of buses (black), branches (gray), and cycles (dashed), while the lower depicts $\errinj$ (black) and $\errflow$ (gray).}
\label{cluster_DE}
\vspace{-0.5cm}
\end{figure}
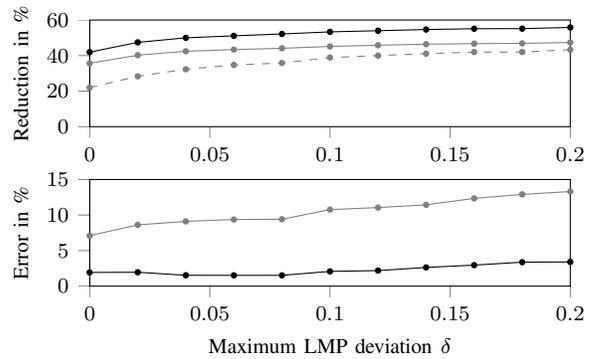
Besides this trade-off, we observed that the error caused by large values of $\tau$ is usually due to a comparably small number of generators which we will refer to as \textit{critical generators}. This indicates that the reduced model does not represent the environment around these generators sufficiently accurate.
As mentioned earlier, the notion of features provides the means to include this knowledge and improve the accuracy of the reduction: The environment of the critical generators is added to the set of features and the reduction is repeated.

This \textit{feature refinement} requires two parameters: an absolute limit of the dispatch error in MW to define when a generator is considered critical (here fixed as 10 MW) and the depth $\vartheta \in \mathbb{N}$. All buses which can be reached from critical generators by traversing a maximum of $\vartheta$ branches are added to the set of features. As shown in Fig.~\ref{1shot_DE}, the error decreases significantly faster than the achieved reduction. This feature refinement allows for large values of $\tau$ with a considerable reduction potential while retaining small dispatch and branch flow errors.

The \textit{market-based reduction} offers the parameter $\delta > 0$. As with the electrical coupling-based reduction, it provides a trade-off between the reduction and the error: A higher $\delta$ enables a more extensive reduction but typically invites a higher error. The results for different $\delta$ are shown in Fig~\ref{cluster_DE}.

\begin{figure}[!t]
\centering
\tikzsetnextfilename{DE_siec}
\setlength\figurewidth{0.9\columnwidth}
\begin{tikzpicture}

\begin{groupplot}[group style={vertical sep=20pt, group size=1 by 2}]
\nextgroupplot[
width=\figurewidth,
height=\figureheight,
label style = {font=\footnotesize}, 
ylabel={Reduction in \%},
xmin=1, xmax=5,
ymin=0, ymax=60,
ticklabel style = {font=\footnotesize},
xtick={2,3,4,5},
xticklabels={(i),(ii.a),(ii.b),(iii)},
tick align=outside,
tick pos=left,
x grid style={white!69.01960784313725!black},
y grid style={white!69.01960784313725!black}
]
\addplot [black, mark=*, mark size=1, mark options={solid}]
table {%
1 0
2 13.3858267716535
3 49.0157480314961
4 41.9291338582677
5 52.1653543307087
};
\addplot [white!50!black, mark=*, mark size=1, mark options={solid}]
table {%
1 0
2 10.8242753623188
3 42.2554347826087
4 35.6884057971014
5 44.1123188405797
};
\addplot [white!50!black, dashed, mark=*, mark size=1, mark options={solid}]
table {%
1 0
2 0.21978021978022
3 30.5494505494506
4 21.978021978022
5 35.8241758241758
};
\path [draw=black, fill opacity=0] (axis cs:1,0)
--(axis cs:1,60);

\nextgroupplot[
width=\figurewidth,
height=\figureheight,
label style = {font=\footnotesize}, 
xlabel={Reduction step},
ylabel={Error in \%},
xmin=1, xmax=5,
ymin=0, ymax=20,
ticklabel style = {font=\footnotesize},
xtick={2,3,4, 5},
xticklabels={(i),(ii.a),(ii.b),(iii)},
tick align=outside,
tick pos=left,
x grid style={white!69.01960784313725!black},
y grid style={white!69.01960784313725!black}
]
\addplot [black, mark=*, mark size=1, mark options={solid}]
table {%
1 0
2 0.118546453899469
3 8.71225273114225
4 1.92266704062884
5 1.51848142740992
};
\addplot [white!50!black, mark=*, mark size=1, mark options={solid}]
table {%
1 0
2 0.070008897005993
3 12.4947782079961
4 7.10874368838231
5 9.4062221797721
};
\path [draw=black, fill opacity=0] (axis cs:1,0)
--(axis cs:1,30);

\end{groupplot}

\end{tikzpicture}
\caption{Error evolution over the reduction steps from Section~\ref{combined}. The upper figure shows the reduction of buses (black), branches (gray), and cycles (dashed), while the lower depicts $\errinj$ (black) and $\errflow$ (gray).}
\vspace{-0.6cm}
\label{de_siec}
\end{figure}
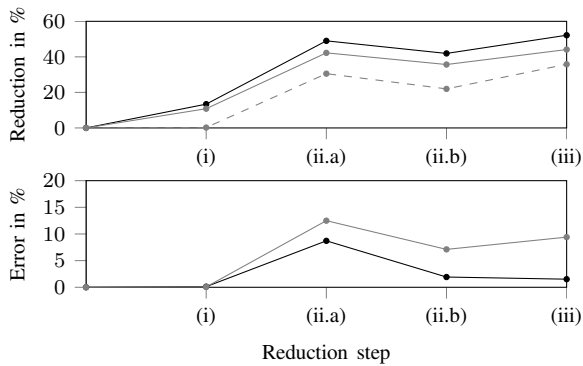

To achieve the targeted dispatch error of $\errinj < 2\%$, the parameters are set to ${\tau = 0.05}, \vartheta = 4$,  and $\delta = 0.08$ which leads to a reduced model with 729 buses ($52.1\%$ reduced), 1234 branches ($44.1\%$ reduced) and 292 cycles\footnote{The number of cycles is relevant in certain expansion strategies, e.g. the \emph{hybrid architecture}~\cite{Hotz:2016ab,Hotz:2018aa,Hotz:2017aa}.} ($35.8\%$ reduced). The error evolution over the reduction steps is shown in Fig.~\ref{de_siec} and results in $\errinj = 1.5\%$ and $\errflow = 9.4\%$.

\section{Application Examples}
\label{sec:verification}
\vspace{-1mm}
In this section, the proposed reduction strategy is applied to three different transmission grids, targeting a reduced model with approximately half the number of buses. This is an exemplary choice -- for a higher accuracy, a lower reduction may be targeted. Furthermore, it is illustrated how the reduction accuracy under the reference scenario for the reduction decisions transfers to different load scenarios. In the following, the German grid (380 and 220\,kV, 455 cycles) of the previous section as well as the Polish and French grid are considered. The Polish grid model represents the Polish 400, 220 and 110\,kV networks during winter peak conditions (2383 buses, 2896 branches and 504 cycles)~\cite{Josz:2018aa}. The French grid model represents the French 380, 225 and 150\,kV networks (2848 buses, 3776 branches and 595 cycles)~\cite{Korab:2018aa}. 
\vspace{-0.2cm}
\begin{table}[!t]
\caption{Reduction and Results}
\begin{center}
\renewcommand{\arraystretch}{\tabstretch}
\setlength{\tabcolsep}{0.85em}
\begin{tabular}{|c|c|c|c|}
\hline
 & {\textbf{German}} & {\textbf{Polish}$^{\mathrm{a}}$} & {\textbf{French}$^{\mathrm{b}}$} \\
\hline
Bus reduction &52.1\,\si{\%} &57.6\,\si{\%} &57.9\,\si{\%}  \\
Branch reduction &44.1\,\si{\%} &48.7\,\si{\%} &50.0\,\si{\%} \\
Cycle reduction &35.8\,\si{\%} &16.2\,\si{\%} &32.9\,\si{\%}\\
\hline
$\errinj$ &\,\,1.5\,\si{\%} &\,\,1.0\,\si{\%} &\,\,0.4\,\si{\%} \\
$\errflow$ &\,\,9.4\,\si{\%} &\,\,8.9\,\si{\%} &\,\,7.8\,\si{\%} \\
\hline
\multicolumn{4}{l}{$^{\mathrm{a}}$ Parameters: $\tau = 0.05, \vartheta = 4$ and $\delta = 0.8$.}\\
\multicolumn{4}{l}{$^{\mathrm{b}}$Parameters: $\tau = 0.03, \vartheta = 7$ and $\delta = 0.1$.}
\end{tabular}
\label{tab1}
\end{center}
\vspace{-0.5cm}
\end{table}

\begin{figure}[!t]

\centering
\tikzsetnextfilename{german_scenarios}
\setlength\figurewidth{0.9\columnwidth}
\begin{tikzpicture}

\definecolor{color0}{rgb}{0.0156862745098039,0.454901960784314,0.584313725490196}

\begin{axis}[
width=\figurewidth,
height=\figureheight,
label style = {font=\footnotesize}, 
ylabel={Error in \%},
xmin=0, xmax=96,
ymin=0, ymax=25,
tick align=outside,
tick pos=left,
ticklabel style = {font=\footnotesize}, 
ytick={5,  15, 25},
yticklabels={$5$, $15$, $25$},
x grid style={white!69.01960784313725!black},
y grid style={white!69.01960784313725!black}
]ll align={left}
]

\addplot [black, mark=*, mark size=0.3, mark options={solid}]
table {%
0 1.45148868092208
1 1.84793004272545
2 3.15984471354452
3 4.52533948039737
4 5.04010933095483
5 5.04010933086471
6 4.52533948038406
7 1.78321266656287
8 4.23900131426052
9 3.08393144926414
10 3.06961288984233
11 3.06961288988399
12 3.08393144927164
13 3.08393144926508
14 3.08393144926838
15 2.86566528202944
16 3.51931621464943
17 2.36389397792979
18 1.45148868082342
19 1.45148868087118
20 3.06961288988636
21 3.72007295024422
22 3.44199791041977
23 1.60398549225547
24 3.15984471360948
25 2.90569563438603
26 1.54574365723273
27 1.66553185540396
28 2.10915124567798
29 2.7698567250116
30 2.41762015094756
31 2.10915124559572
32 1.52118681781951
33 5.15100815738659
34 3.17576998385271
35 3.71843609161475
36 3.72007295024535
37 3.71843609160786
38 3.175769983839
39 3.49984175471118
40 3.49984175468812
41 3.7200729502175
42 1.45148868083288
43 2.36389397561338
44 3.06655418367777
45 3.51931621463872
46 3.2740376606554
47 3.49984175468812
48 3.38654036573478
49 2.76985672495137
50 4.5253394803259
51 5.55293986868033
52 6.49888802921579
53 6.49888803030834
54 5.55293986853267
55 2.76985672495137
56 5.00875453991878
57 3.49984175468812
58 3.08393144926446
59 2.36389397561441
60 1.45148868082868
61 2.36389396619623
62 1.45148868082926
63 1.45148868082868
64 3.06655418367555
65 3.06961288986649
66 3.06961288986649
67 2.8656652820167
68 3.27403766091616
69 3.27403766091716
70 2.86566528202115
71 3.49984175469974
72 1.54574365722898
73 1.78321266656652
74 1.52118681782233
75 2.10915124567798
76 2.41762015095333
77 2.76985672495805
78 3.58022387257983
79 3.58022387257983
80 2.10915124559572
81 3.38654036610659
82 4.23900131427165
83 3.72007295024535
84 2.86566528202944
85 2.86566528202115
86 3.27403766091716
87 3.72007295018625
88 3.72007295018664
89 3.27403766091598
90 3.51931621464032
91 3.08393144926414
92 3.08393144926414
93 1.45148868087334
94 2.86566528202195
95 3.17576998385455
96 5.15100815732939
};
\addplot [white!69.01960784313725!black, mark=*, mark size=0.3, mark options={solid}]
table {%
0 9.39864122281347
1 9.57612911188153
2 11.6074750023168
3 14.4524863863763
4 15.7325670537347
5 15.7325670537203
6 14.4524863863776
7 9.64196700859938
8 13.2315105950682
9 12.3140140377763
10 12.3532567445313
11 12.3532567445808
12 12.3140140377745
13 12.3140140377685
14 12.3140140377745
15 11.8515100738452
16 12.5649284131742
17 11.0178329426908
18 9.39864122252774
19 9.39864122268513
20 12.3532567445825
21 12.4978929886017
22 12.5229044871168
23 9.40160623788114
24 11.607475002389
25 11.8706881378805
26 9.27847550817058
27 9.35221236108799
28 9.93799173747212
29 10.9409807049482
30 10.4134289806443
31 9.93799173744002
32 9.2293567745911
33 12.3056170298036
34 12.5088226818188
35 12.4213669884067
36 12.4978929886027
37 12.4213669884052
38 12.5088226818062
39 13.3213903055292
40 13.321390305503
41 12.4978929885611
42 9.39864122269984
43 11.0178329422357
44 12.0581507992813
45 12.5649284131606
46 12.1849074696831
47 13.321390305503
48 12.4012627259336
49 10.9409807047858
50 14.4524863864637
51 17.0633058623806
52 19.3363937427118
53 19.3363937461076
54 17.0633058621436
55 10.9409807047858
56 12.5806297047251
57 13.321390305503
58 12.3140140377728
59 11.0178329422331
60 9.39864122269472
61 11.0178329413335
62 9.39864122269378
63 9.39864122269472
64 12.0581507992794
65 12.3532567445617
66 12.3532567445617
67 11.8515100738309
68 12.1849074698018
69 12.1849074698043
70 11.8515100738389
71 13.321390305517
72 9.27847550817063
73 9.64196700859836
74 9.22935677460483
75 9.93799173747212
76 10.4134289806925
77 10.9409807048284
78 12.4092063376024
79 12.4092063376024
80 9.93799173744002
81 12.4012627260447
82 13.2315105950923
83 12.4978929886027
84 11.8515100738452
85 11.8515100738389
86 12.1849074698043
87 12.4978929884995
88 12.4978929885009
89 12.1849074698046
90 12.5649284131735
91 12.3140140377763
92 12.3140140377763
93 9.39864122267526
94 11.8515100738433
95 12.5088226818587
96 12.3056170297901
};

\path [draw=black, fill opacity=0] (axis cs:0,0)
--(axis cs:96,0);

\end{axis}

\end{tikzpicture}%
\tikzsetnextfilename{polish_scenarios}
\begin{tikzpicture}

\begin{axis}[
width=\figurewidth,
height=\figureheight,
label style = {font=\footnotesize}, 
ylabel={Error in \%},
xmin=0, xmax=96,
ymin=0, ymax=25,
tick align=outside,
tick pos=left,
ticklabel style = {font=\footnotesize}, 
ytick={5,  15, 25},
yticklabels={$5$, $15$, $25$},
x grid style={white!69.01960784313725!black},
y grid style={white!69.01960784313725!black}
]
\addplot [black, mark=*, mark size=0.3, mark options={solid}]
table {%
0 1.05675146104105
1 0.201503750702613
2 0.19419854831857
3 0.221445586074455
4 0.214989313484942
5 0.214989313485291
6 0.221445586071939
7 0.563391601075956
8 1.71403947257546
9 0.967916548282671
10 1.21049462564391
11 1.21049462564544
12 0.967916548287577
13 0.967916548282671
14 0.967916548281591
15 1.2757328797872
16 0.91898300573439
17 1.00847443154988
18 1.05675146104105
19 1.05675146104646
20 1.21049462563798
21 1.322580792664
22 0.53988221769616
23 0.322134117980272
24 0.194198548317407
25 0.220782925467927
26 0.285069006991351
27 0.21017096603559
28 0.198865865175188
29 0.194783821609419
30 0.198590488523034
31 0.19886586517586
32 0.220721603387323
33 0.351568298552167
34 1.70553217474955
35 1.34082891553999
36 1.322580792664
37 1.34082891559998
38 1.70553217474312
39 1.74291190122439
40 1.74291190120958
41 1.32258079268904
42 1.05675146104762
43 1.00847443155045
44 1.15534892780564
45 0.918983005715526
46 1.34188958092012
47 1.74291190122066
48 0.432391576678306
49 0.194783821610525
50 0.221445586075652
51 0.36339753193375
52 0.941960925940095
53 0.941960925890055
54 0.36339753193375
55 0.194783821610671
56 0.473117348069807
57 1.7429119012133
58 0.967916548282671
59 1.00847443155067
60 1.05675146104222
61 1.00847443155045
62 1.05675146104105
63 1.05675146104105
64 1.15534892780564
65 1.21049462564544
66 1.21049462564544
67 1.27573287981093
68 1.34188958092554
69 1.34188958092554
70 1.2757328797872
71 1.74291190122066
72 0.285069006991351
73 0.563391601054614
74 0.220721603387323
75 0.19886586517586
76 0.198590488523073
77 0.194783821609419
78 0.200944465869135
79 0.200944465871548
80 0.198865865175032
81 0.432391576678306
82 1.71403947257546
83 1.322580792664
84 1.27573287980855
85 1.27573287981093
86 1.34188958093245
87 1.32258079268904
88 1.32258079268904
89 1.34188958092554
90 0.918983005715526
91 0.967916548282671
92 0.967916548281591
93 1.05675146104762
94 1.2757328797872
95 1.70553217473732
96 0.351568298548677
};
\addplot [white!69.01960784313725!black, mark=*, mark size=0.3, mark options={solid}]
table {%
0 8.97990224221517
1 11.6240770849436
2 12.1869532474398
3 12.7750116065009
4 13.0070227425708
5 13.0070227425727
6 12.7750116064496
7 11.6261362061129
8 11.2144041705009
9 8.96798485667586
10 9.05023161356291
11 9.05023161356586
12 8.96798485667985
13 8.96798485667586
14 8.9679848566749
15 9.42020315857532
16 9.04229765223248
17 9.01094750022071
18 8.97990224221517
19 8.97990224221919
20 9.05023161355052
21 9.83607846707236
22 10.5562840795438
23 11.7548915880402
24 12.1869532474404
25 10.6516611942128
26 11.9215705829815
27 11.6894022734235
28 11.6747760989725
29 12.0712556257267
30 11.74841555793
31 11.674776098973
32 11.5533743210953
33 10.6250694714455
34 11.0410665591425
35 10.1632776866348
36 9.83607846707236
37 10.1632776866743
38 11.0410665591235
39 11.1603287786177
40 11.160328778589
41 9.83607846708454
42 8.97990224222339
43 9.01094750022709
44 8.98399662276951
45 9.04229765222686
46 9.6325493205549
47 11.1603287785971
48 10.6362986106846
49 12.071255625728
50 12.7750116065039
51 13.4743184940553
52 14.7332589074918
53 14.7332589073098
54 13.4743184940553
55 12.0712556257264
56 10.9565195636973
57 11.1603287786096
58 8.96798485667586
59 9.01094750022985
60 8.97990224221938
61 9.01094750022709
62 8.97990224221517
63 8.97990224221517
64 8.98399662276951
65 9.05023161356586
66 9.05023161356586
67 9.4202031585982
68 9.63254932056597
69 9.63254932056597
70 9.42020315857532
71 11.1603287785971
72 11.9215705829815
73 11.6261362058129
74 11.5533743210953
75 11.674776098973
76 11.7484155579305
77 12.0712556257267
78 12.3353320450211
79 12.335332045005
80 11.6747760989724
81 10.6362986106846
82 11.2144041705009
83 9.83607846707236
84 9.42020315859638
85 9.4202031585982
86 9.63254932056918
87 9.83607846708454
88 9.83607846708454
89 9.63254932056597
90 9.04229765222686
91 8.96798485667586
92 8.9679848566749
93 8.97990224222339
94 9.42020315857532
95 11.0410665590995
96 10.6250694714386
};
\path [draw=black, fill opacity=0] (axis cs:0,0)
--(axis cs:96,0);

\end{axis}

\end{tikzpicture}%
\tikzsetnextfilename{french_scenarios}
\begin{tikzpicture}

\definecolor{color0}{rgb}{0.0156862745098039,0.454901960784314,0.584313725490196}

\begin{axis}[
width=\figurewidth,
height=\figureheight,
label style = {font=\footnotesize}, 
xlabel={Scenario},
ylabel={Error in \%},
xmin=0, xmax=96,
ymin=0, ymax=25,
tick align=outside,
tick pos=left,
ticklabel style = {font=\footnotesize}, 
ytick={5,  15, 25},
yticklabels={$5$, $15$, $25$},
x grid style={white!69.01960784313725!black},
y grid style={white!69.01960784313725!black}
]

\addplot [black, mark=*, mark size=0.3, mark options={solid}]
table {%
0 0.466254304501672
1 0.363122481241007
2 0.363435704995805
3 0.261455172245956
4 0.25971547105158
5 0.25971547105158
6 0.261455172249071
7 0.435095893250428
8 0.68305343990433
9 0.804792639344634
10 0.503396336250782
11 0.503396336251795
12 0.804792639349781
13 0.804792639343784
14 0.804792639349781
15 0.515920941688469
16 0.576965628235389
17 0.465964810659204
18 0.466254304502897
19 0.466254304501982
20 0.503396336250782
21 0.486202794973142
22 0.659114483550489
23 0.410978456115062
24 0.363435704995497
25 0.448264358879772
26 0.409706205048738
27 0.358428566857374
28 0.336965484524132
29 0.355426246655139
30 0.341977270064025
31 0.336965484520777
32 0.390331521109273
33 0.455849567789104
34 0.498334721307672
35 0.486197593786059
36 0.486202794976805
37 0.486197593785787
38 0.498334721307672
39 0.462394755899587
40 0.462394755899587
41 0.486202794976805
42 0.466254304502897
43 0.465964810659204
44 0.472685270155447
45 0.576965628235389
46 0.516555069163169
47 0.462394755899589
48 0.459342200817577
49 0.355426246655139
50 0.261455172251964
51 0.259757272088457
52 0.248311645362842
53 0.248311645362842
54 0.259757272088457
55 0.355426246655139
56 0.458695818686751
57 0.462394755899093
58 0.804792639349781
59 0.465964810325454
60 0.466254304502977
61 0.465964810322842
62 0.466254304502897
63 0.466254304502897
64 0.472685270155411
65 0.503396336249268
66 0.503396336250782
67 0.515920941688375
68 0.51655506916066
69 0.516555069163169
70 0.515920941688375
71 0.462394755899093
72 0.409706205048214
73 0.435095893245372
74 0.390331521111591
75 0.336965484520777
76 0.341977270061945
77 0.355426246655139
78 0.3666679179514
79 0.3666679179514
80 0.336965484519955
81 0.459342200802333
82 0.683053439971429
83 0.486202794977237
84 0.515920941687569
85 0.515920941687569
86 0.51655506915975
87 0.486202794976805
88 0.486202794977237
89 0.51655506915975
90 0.576965628235389
91 0.804792639349781
92 0.804792639349781
93 0.466254304501903
94 0.515920941688469
95 0.498334721307161
96 0.455849567788402
};
\addplot [white!69.01960784313725!black, mark=*, mark size=0.3, mark options={solid}]
table {%
0 7.84310702730508
1 8.23166284620149
2 8.03050785275682
3 7.46013418059305
4 7.74836334742513
5 7.74836334742513
6 7.46013418059673
7 7.7955784000923
8 8.11761116418741
9 8.30101459688916
10 7.84350875382687
11 7.84350875382786
12 8.30101459688624
13 8.30101459688341
14 8.30101459688624
15 7.85494623313208
16 7.95419197227787
17 7.64642065984463
18 7.84310702730579
19 7.84310702730301
20 7.84350875382687
21 7.71046498627035
22 7.80978367441982
23 7.97741693438683
24 8.03050785275724
25 7.72234329310646
26 8.01882277426262
27 7.90787978186015
28 8.32599973240305
29 8.50531041394437
30 8.58780172569925
31 8.3259997324017
32 8.03901634202217
33 7.4198214122923
34 8.2415335351452
35 7.900582725528
36 7.71046498627261
37 7.90058272552939
38 8.2415335351452
39 7.9542452493281
40 7.9542452493281
41 7.71046498627261
42 7.84310702730579
43 7.64642065984463
44 7.92862017993504
45 7.95419197227787
46 7.91027265795018
47 7.9542452493257
48 7.62891445358547
49 8.50531041394437
50 7.46013418060151
51 8.12815321200722
52 8.09099456432157
53 8.09099456432157
54 8.12815321200722
55 8.50531041394437
56 7.6026110758084
57 7.95424524932584
58 8.30101459688624
59 7.64642065902881
60 7.84310702730605
61 7.64642065902934
62 7.84310702730579
63 7.84310702730579
64 7.92862017993582
65 7.84350875382706
66 7.84350875382687
67 7.85494623312562
68 7.91027265794579
69 7.91027265795018
70 7.85494623312562
71 7.95424524932584
72 8.01882277426069
73 7.79557840009126
74 8.03901634202911
75 8.3259997324017
76 8.58780172569674
77 8.50531041394437
78 7.88244096262746
79 7.88244096262746
80 8.32599973240063
81 7.6289144536319
82 8.11761116420755
83 7.71046498627247
84 7.85494623312585
85 7.85494623312585
86 7.91027265794574
87 7.71046498627261
88 7.71046498627247
89 7.91027265794574
90 7.95419197227787
91 8.30101459688624
92 8.30101459688624
93 7.84310702730275
94 7.85494623313208
95 8.24153353514474
96 7.41982141229164
};
\path [draw=black, fill opacity=0] (axis cs:0,0)
--(axis cs:96,0);

\end{axis}

\end{tikzpicture}%
\caption{Dispatch error $\errinj$ (black) and branch flow error $\errflow$ (gray) of the reduced model of the German (top), Polish (middle), and French (bottom) transmission grid.}
\vspace{-0.5cm}
\label{scenarios}
\end{figure}
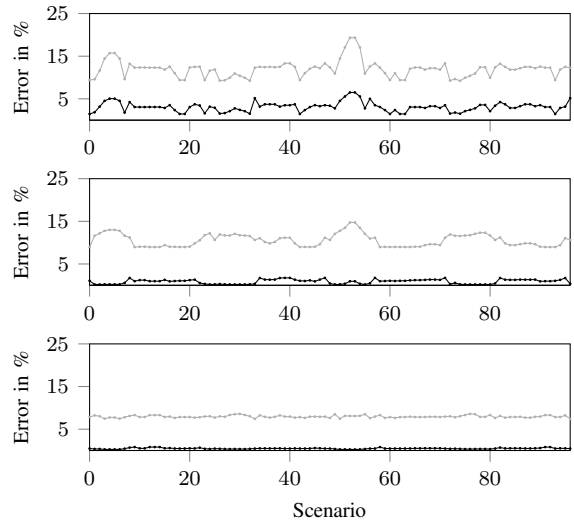
\subsection{Reduction and Results for the Reference Scenario}\label{AA} 
As documented in Table~\ref{tab1}, reduced models with a similar reduction extent and accuracy to that of the German grid discussed in Section~\ref{combined} can be constructed for the Polish and French system. This is a consequence of the fact that the reduction strategy can be adapted to different grids by adjusting the reduction parameters. On account of the different approaches, there is an inherent adaptability: if one reduction approach is less effective for a certain grid, the other approaches can be applied more radically to achieve a similar reduction.

\subsection{Verification for Different Load Scenarios}\label{AB}
The reduction strategy is solely based on a reference scenario of the grid. However, for conversion planning it is necessary that the reduced model offers also sufficient accuracy under different load scenarios. In the following, this is examined for all three grids. The different load scenarios are generated by scaling the loads of the reference scenario according to the 96 hours of the exemplary winter and summer weekday and weekend presented by the IEEE Reliability Test System Task Force~\cite[Table~4]{Grigg:1999aa}.

As shown in Fig.~\ref{scenarios}, the generalization is particular to a model: the French grid, which is the most strongly meshed among these examples, exhibits an especially consistent behavior. The German and Polish grid fluctuate somewhat more and, depending on the application, a slightly more conservative reduction may be necessary. Still, all reduced models exhibit an adequate accuracy for all scenarios, which can be attributed to the exclusive reduction of subgrids with a limited impact on the overall system behavior.

\section{Conclusion}
\label{sec:conclusion}

The proposed reduction strategy combines a topology-based, electrical coupling-based, and market-based subgrid selection with a specifically designed reduction method to produce a reduced grid model which preserves designated features and the structure of the original grid. This approach avoids the addition of artificial entities and the implementation of physically invalid or application-adverse reduction measures by being aware of and preserving features and the relation of every entity to its physical counterpart.

Inherently, on account of its adjustable parameters, the reduction strategy offers a trade-off between the reduced model's complexity and accuracy, enabling the construction of a reduced model with a specific accuracy. The combination of the different approaches also allows the strategy to adapt to different grids: if one reduction approach proves to be less effective for a certain grid, the targeted reduction can be achieved by applying the others more drastically. The reduction strategy offers considerable reduction potential, e.g. to half of the number of buses at a very moderate error, which enables power flow studies at substantially lower computational costs.

Due to the careful selection of subgrids for the reduction and the preservation of features, the accuracy for the reference scenario generalizes well to different load scenarios: An appropriate parametrization under the reference scenario leads to a reduced model which adequately represents the original grid model for various load scenarios.

Especially in the context of comprehensive power flow studies for capacity expansion it is of the essence that conclusions drawn from the reduced grid model transfer to all load cases in order to arrive at reasonable capacity expansion measures. Moreover, the objective of these studies is to identify areas of the grid which require capacity expansion. The proposed reduction strategy is especially qualified in this case as it is designed to preserve such structures within the reduced model and thereby draws the focus to areas that are crucial to the overall grid behavior.
\vspace{-0.23cm}

\bibliographystyle{IEEEtran}
\bibliography{literature}

\begin{thebibliography}{10}
\providecommand{\url}[1]{#1}
\csname url@samestyle\endcsname
\providecommand{\newblock}{\relax}
\providecommand{\bibinfo}[2]{#2}
\providecommand{\BIBentrySTDinterwordspacing}{\spaceskip=0pt\relax}
\providecommand{\BIBentryALTinterwordstretchfactor}{4}
\providecommand{\BIBentryALTinterwordspacing}{\spaceskip=\fontdimen2\font plus
\BIBentryALTinterwordstretchfactor\fontdimen3\font minus
  \fontdimen4\font\relax}
\providecommand{\BIBforeignlanguage}[2]{{%
\expandafter\ifx\csname l@#1\endcsname\relax
\typeout{** WARNING: IEEEtran.bst: No hyphenation pattern has been}%
\typeout{** loaded for the language `#1'. Using the pattern for}%
\typeout{** the default language instead.}%
\else
\language=\csname l@#1\endcsname
\fi
#2}}
\providecommand{\BIBdecl}{\relax}
\BIBdecl

\bibitem{Capros:2013aa}
P.~P. Capros, A.~D. Vita, N.~Tasios, D.~Papadopoulos, P.~Siskos, E.~Apostolaki
  \emph{et~al.}, ``{EU} {E}nergy, {T}ransport and {GHG} {E}missions -- {T}rends
  to 2050,'' European Commission, Tech. Rep., Dec. 2013.

\bibitem{Flechtner:2015aa}
J.~Flechtner and D.~S. Bolay, ``Faktenpapier {A}usbau der {S}tromnetze,'' DIHK
  - Deutscher Industrie- und Handelskammertag Berlin/Br{\"u}ssel, Tech. Rep.,
  Jan. 2015.

\bibitem{Netzentwicklungsplan2017a}
{50Hertz Transmission}, {Amprion}, {TenneT TSO}, and {TransnetBW},
  ``{N}etzentwicklungsplan {S}trom 2030 (2.~{E}ntwurf),'' Tech. Rep., May 2017.

\bibitem{Papaemmanouil:2011aa}
A.~Papaemmanouil and G.~Andersson, ``On the reduction of large power system
  models for power market simulations,'' in \emph{Proc. 17th Power Systems
  Computation Conference (PSCC)}, Aug. 2011.

\bibitem{Liacco:1978aa}
T.~E.~D. Liacco, S.~C. Savulescu, and K.~A. Ramarao, ``An on-line topological
  equivalent of a power system,'' \emph{IEEE Trans. on Power Apparatus and
  Syst.}, vol. PAS-97, no.~5, pp. 1550--1563, Sep. 1978.

\bibitem{Bergen:2000aa}
A.~R. Bergen and V.~Vittal, \emph{Power Systems Analysis}, 2nd~ed.\hskip 1em
  plus 0.5em minus 0.4em\relax Prentice Hall Inc., 2000.

\bibitem{Savulescu:2002aa}
S.~C. Savulescu, ``Solving open access transmission and security analysis
  problems with the short-circuit currents method,'' in \emph{Pennwell Latin
  America Power Conf., Controlling and Automation Energy Session}, Aug. 2002.

\bibitem{Ward:1949aa}
J.~B. Ward, ``Equivalent circuits for power-flow studies,'' \emph{AIEE
  Transactions}, vol.~68, pp. 373--382, Jul. 1949.

\bibitem{Wu:1983aa}
F.~Wu and A.~Monticelli, ``Critical review of external network modelling for
  online security analysis,'' \emph{Int. J. Electrical Power \& Energy Syst.},
  vol.~5, no.~4, 1983.

\bibitem{Brown:1985aa}
H.~E. Brown, \emph{Solution of Large Networks by Matrix Methods}.\hskip 1em
  plus 0.5em minus 0.4em\relax New York: John Wiley \& Sons Inc, 1985.

\bibitem{Duran:1972aa}
H.~Duran and N.~Arvanitidis, ``Simplification for area security analysis: {A}
  new look at equivalencing,'' \emph{IEEE Trans. Power App. Syst.}, vol.
  PAS-91, pp. 670--679, Mar. 1972.

\bibitem{Monticelli:1979aa}
A.~Monticelli, S.~Deckmann, A.~Garcia, and B.~Stott, ``Real-time external
  equivalents for static security analysis,'' \emph{IEEE Trans. Power App.
  Syst.}, vol. PAS-98, no.~2, pp. 498--508, Mar. 1979.

\bibitem{Deckmann:1980aa}
S.~Deckmann, A.~Pizzolante, A.~Monticelli, B.~Stott, and O.~Alsac, ``Studies on
  power system load flow equivalencing,'' \emph{IEEE Trans. Power App. Syst.},
  vol. PAS-99, no.~6, pp. 2301--2310, Nov. 1980.

\bibitem{Dimo:1975aa}
P.~Dimo, \emph{Nodal Analysis of Power Systems}.\hskip 1em plus 0.5em minus
  0.4em\relax Kent, England: Abacus Press, 1975.

\bibitem{Gavrilas:2008aa}
M.~Gavrilas, O.~Ivanov, and G.~Gavrilas, ``{REI} equivalent design for electric
  power systems with genetic algorithms,'' \emph{WSEAS Trans. Circuits and
  Syst.}, vol.~7, no.~10, pp. 911--921, Oct. 2008.

\bibitem{Milano:2009aa}
F.~Milano and K.~Srivastava, ``Dynamic {REI} equivalents for short circuit and
  transient stability analyses,'' \emph{Electric Power Systems Research},
  vol.~79, no.~6, pp. 878 -- 887, Jan. 2009.

\bibitem{Tinney:1983aa}
W.~F. Tinney and W.~L. Powell, ``The {REI} approach to power network
  equivalents,'' in \emph{Proc. IEEE PICA Conf.}, Toronto, Canada, May 1977,
  pp. 314--320.

\bibitem{Housos:1980aa}
E.~C. Housos, G.~Irisarri, R.~M. Porter, and A.~M. Sasson, ``Steady state
  network equivalents for power system planning applications,'' vol. PAS-99,
  no.~6, pp. 2113--2120, Nov. 1980.

\bibitem{Singh:2005aa}
H.~K. Singh and S.~C. Srivastava, ``A reduced network representation suitable
  for fast nodal price calculations in electricity markets,'' in \emph{Proc.
  IEEE Power Engineering Society General Meeting}, vol.~2, Jun. 2005, pp.
  2070--2077.

\bibitem{Shi:2015aa}
D.~Shi and D.~J. Tylavsky, ``A novel bus-aggregation-based structure-preserving
  power system equivalent,'' \emph{IEEE Trans. Power Syst.}, vol.~30, no.~4,
  pp. 1977--1986, Jul. 2015.

\bibitem{Hotz:2019aa}
\BIBentryALTinterwordspacing
M.~Hotz \emph{et~al.}, ``{\emph{hynet:}} {A}n optimal power flow framework for
  hybrid {AC}/{DC} power systems (v1.1.0),'' Mar. 2019. [Online]. Available:
  \url{http://gitlab.com/tum-msv/hynet}
\BIBentrySTDinterwordspacing

\bibitem{Hotz:2018ab}
\BIBentryALTinterwordspacing
M.~Hotz and W.~Utschick, ``\textit{{hynet}:} {A}n optimal power flow framework
  for hybrid {AC}/{DC} power systems,'' \emph{arXiv:1811.10496}, Nov. 2018.
  [Online]. Available: \url{http://arxiv.org/abs/1811.10496v1}
\BIBentrySTDinterwordspacing

\bibitem{Hotz:2018aa}
------, ``The hybrid transmission grid architecture: {B}enefits in nodal
  pricing,'' \emph{IEEE Trans. Power Syst.}, vol.~33, no.~2, Mar. 2018.

\bibitem{Hewes:2016aa}
D.~Hewes, S.~Altschaeffl, I.~Boiarchuk, and R.~Witzmann, ``Development of a
  dynamic model of the {E}uropean transmission system using publicly available
  data,'' in \emph{Proc. IEEE Int. Energy Conf. (ENERGYCON)}, Apr. 2016.

\bibitem{Hotz:2016ab}
M.~Hotz and W.~Utschick, ``A hybrid transmission grid architecture enabling
  efficient optimal power flow,'' \emph{IEEE Trans. Power Syst.}, vol.~31,
  no.~6, pp. 4504--4516, Nov. 2016.

\bibitem{Hotz:2017aa}
M.~Hotz, W.~Utschick, D.~Hewes, I.~Boiarchuk, and R.~Witzmann, ``Reducing the
  need for new lines in {G}ermany's energy transition: {T}he hybrid
  transmission grid architecture,'' \emph{Int. ETG Congress 2017}, Nov. 2017.

\bibitem{Josz:2018aa}
C.~Josz, S.~Fliscounakis, J.~Maeght, and P.~Panciatici, ``{AC} {P}ower {F}low
  {D}ata in {MATPOWER} and {QCQP} format: i{T}esla, {RTE} {S}napshots, and
  {PEGASE},'' Tech. Rep., Jun. 2018.

\bibitem{Korab:2018aa}
R.~Korab, ``{IEEE} {PES} {P}ower {G}rid {L}ibrary - case2383wp,'' Tech. Rep.,
  Jun. 2018.

\bibitem{Grigg:1999aa}
C.~Grigg, P.~Wong, P.~Albrecht, R.~Allan \emph{et~al.}, ``The {IEEE}
  reliability test system - 1996. {A} report prepared by the {R}eliability
  {T}est {S}ystem {T}ask {F}orce of the {A}pplication of {P}robability
  {M}ethods {S}ubcommittee,'' \emph{IEEE Trans. Power Syst.}, vol.~14, no.~3,
  pp. 1010--1020, Aug. 1999.

\end{thebibliography}

\end{document}